\documentclass[11pt, reqno]{amsart}
\usepackage{amsmath, amsthm, amscd, amsfonts, amssymb, graphicx, color}
\usepackage[bookmarksnumbered, colorlinks, plainpages]{hyperref}
\hypersetup{colorlinks=true,linkcolor=red, anchorcolor=green, citecolor=cyan, urlcolor=red, filecolor=magenta, pdftoolbar=true}

\newtheorem{theorem}{Theorem}[section]
\newtheorem{lemma}[theorem]{Lemma}
\newtheorem{proposition}[theorem]{Proposition}
\newtheorem{corollary}[theorem]{Corollary}
\theoremstyle{definition}

\theoremstyle{remark}
\newtheorem{remark}[theorem]{Remark}
\numberwithin{equation}{section}
\newcommand{\Mn}{\mathbb M_n}

\newcommand{\spn}{\operatorname{span}}
\newcommand{\tr}{\operatorname{Tr}}
\newcommand{\dg}{\operatorname{diag}}
\begin{document}
\title[Unbounded Gaps in Matrix Orders]{Unbounded Gaps Between Matrix Monotonicity and Frobenius Radial Convexity}

\author[M. Kian]{Mohsen Kian}
\address{M. Kian: Department of Mathematics, University of Bojnord, P. O. Box
1339, Bojnord 94531, Iran}
 \email{kian@ub.ac.ir }

\keywords{Matrix monotonicity, Frobenius radial convexity,  Loewner matrices, rational interpolation}
\subjclass[2020]{Primary 47A63; Secondary 15A39, 47A60}

\begin{abstract}
We show that the gap between matrix monotonicity and
Frobenius radial convexity can be arbitrarily large.  For every integer $m\ge 2$, we construct rational functions
that are matrix monotone of order $m$ but not of order $m+1$,
and Frobenius radially convex of order $2m-1$ but not of
order $2m$. Thus, the difference between their exact
matrix-monotonicity and radial-convexity orders is $m-1$,
which can be arbitrarily large. As a consequence, for every
$n\ge 3$, Frobenius radial convexity of order $n$ does not
imply matrix monotonicity of the same order, either
increasing or decreasing.
\end{abstract}

\maketitle

\section{Introduction}
Matrix monotonicity extends ordinary monotonicity to Hermitian matrices. A real function $f:J\subseteq \mathbb{R}\to\mathbb{R}$ is called matrix monotone of order $n$ if it preserves the Loewner order, i.e., if $A\le B$ implies $f(A)\le f(B)$ for all $n\times n$ Hermitian matrices $A$ and $B$  with spectra in $J$.   This requirement becomes more
restrictive as the matrix size increases. Hansen, Ji, and Tomiyama \cite{HJT2004} proved that
these inclusions are strict on every nontrivial interval other than $\mathbb{R}$.
Osaka, Silvestrov, and Tomiyama \cite{OST2007} later constructed polynomials belonging to
the gaps between successive classes.
Conditional positivity of Loewner matrices has also been used to study finite-order matrix monotonicity and convexity \cite{HS2012}, while perturbations of Cauchy matrices have been investigated in \cite{ST2021}.

Czerwi\'nska, Franco, and Merino \cite{CFM2026} studied matrix monotonicity through the Frobenius distance
$$d_f(A,B)=\|f(A)-f(B)\|_{\mathrm F},$$
where $f$ is real-valued on $(0,\infty)$ and $A,B$ are positive definite. A function $f$ is
Frobenius radially convex of order $n$ if
$$A\le C\le B \quad\Longrightarrow\quad d_f(A,C)\le d_f(A,B) \qquad\text{and}\qquad
d_f(C,B)\le d_f(A,B)$$
for positive definite $A,B,C\in\Mn$.
Every matrix-monotone increasing and matrix-monotone decreasing function  of order $n$ is   Frobenius radially convex of order $n$, see \cite[Proposition~3.9]{CFM2026}.

Czerwi\'nska, Franco, and Merino conjectured that Frobenius radial convexity of order $n$
characterizes matrix monotonicity of the same order, in either direction. They proved this
for differentiable functions when $n=2$.

The conjecture was subsequently disproved at order $3$ by a rational counterexample
documented in October 2026 \cite{TRT2026}. This example is the $m=2$ case of the family
studied here, with $\varepsilon=1/20000$.

Let $\mathcal M_n$ and $\mathcal R_n$ denote the classes of matrix-monotone increasing
and Frobenius radially convex functions of order $n$, respectively.
For each $m\ge2$, we show (Theorem~\ref{th-gap})  that the rational function
$f_{m,\varepsilon}$ defined in \eqref{eq--fam} satisfies
$$f_{m,\varepsilon}\in\mathcal M_m\setminus\mathcal M_{m+1}
\qquad\text{and}\qquad f_{m,\varepsilon}\in
\mathcal R_{2m-1}\setminus\mathcal R_{2m}$$
for all sufficiently small $\varepsilon>0$.
The difference between the two exact orders is therefore $m-1$, which is unbounded.
In particular, the conjectured implication fails at every order $n\ge3$, even when
matrix monotonicity is allowed in either direction.

\section{Preliminaries}
Let $\Mn$ denote the algebra of $n\times n$ complex matrices. For Hermitian
matrices $A,B\in\mathbb{M}_n$, we write $A\le B$ if $B-A$ is positive semidefinite,
and $A< B$ if $B-A$ is positive definite.
The identity matrix is denoted by $I$.
The Frobenius norm is given by
\begin{align*}
\|X\|_{\mathrm F}=\bigl(\tr(X^*X)\bigr)^{1/2},\qquad X\in\Mn.
\end{align*}
Let $J\subseteq\mathbb{R}$ be an open interval. A real-valued function $f$ on $J$ is
\emph{matrix monotone of order $n$} if  $A\le B$ implies that $f(A)\le f(B)$
for all Hermitian $A,B\in\Mn$ with spectra in $J$. We denote the class of such functions by
$\mathcal{M}_n(J)$. For $J=(0,\infty)$, we simply write
$\mathcal{M}_n$.

For a differentiable real-valued function $f$ on $J$, its first divided difference is
\begin{align*}
f^{[1]}(s,t)=\begin{cases}
\dfrac{f(s)-f(t)}{s-t}, & s\ne t,\\
f'(s), & s=t.
\end{cases}
\end{align*}
The associated Loewner matrix is
\begin{align*}
L_f(s_1,\ldots,s_n)=\bigl[f^{[1]}(s_i,s_j)\bigr]_{i,j=1}^n.
\end{align*}
By Loewner's criterion (see \cite[Chapter V]{Bha1997}),
a function $f\in C^1(J)$ (continuously differentiable on $J$)
  belongs to $\mathcal{M}_n(J)$ if and only if $L_f(s_1,\ldots,s_n)\geq0$ for every $s_1,\ldots,s_n\in J$.

For a sufficiently smooth function $f$, we denote its divided difference of order
$\ell$ by $f^{[\ell]}(s_0,\ldots,s_\ell)$. These are defined recursively by
\begin{align*}
f^{[\ell]}(s_0,\ldots,s_\ell)=\frac{f^{[\ell-1]}(s_1,\ldots,s_\ell)
-f^{[\ell-1]}(s_0,\ldots,s_{\ell-1})}{s_\ell-s_0},
\end{align*}
when $s_0\ne s_\ell$, with $f^{[0]}=f$. At repeated points, divided differences
are defined by continuity. In particular,
\begin{align*}
f^{[\ell]}(s,\ldots,s)=\frac{f^{(\ell)}(s)}{\ell!}.
\end{align*}

Let $f\in C^1(J)$ and let $C\in\Mn$ be Hermitian with spectrum in $J$.
The Fr\'echet derivative of $f$ at $C$, denoted by $Df(C)$, is the linear map
on Hermitian matrices characterized by
\begin{align*}
f(C+X)=f(C)+Df(C)[X]+o(\|X\|_{\mathrm F})
\end{align*}
as $\|X\|_{\mathrm F}\to0$. By the Daleckii--Krein formula
(see \cite[Chapter V]{Bha1997}),  if $C=U\dg(s_1,\ldots,s_n)U^*$,
then
\begin{align*}
Df(C)[X]=U\bigl(L_f(s_1,\ldots,s_n)\circ(U^*XU)\bigr)U^*
\end{align*}
for every Hermitian $X\in\Mn$, where $\circ$ denotes the Hadamard product.

In particular, $Df(C)$ is self-adjoint with respect to the Hilbert--Schmidt
inner product, i.e.,
\begin{align*}
\tr\bigl(Df(C)[X]Y\bigr)=\tr\bigl(XDf(C)[Y]\bigr)
\end{align*}
for all Hermitian $X,Y\in\Mn$.

For a real function $f$ defined on an open interval $J\subseteq\mathbb R$,
  define
\begin{align*}
d_f(A,B)=\|f(A)-f(B)\|_{\mathrm F}
\end{align*}
for  all  Hermitian matrices $A,B\in\Mn$ with spectra in $J$. Following the general definition of radial convexity in \cite[Definition~1.1]{CFM2026},
we say that $f$ is \emph{Frobenius radially convex
of order $n$ on $J$}
\begin{align}\label{eq--rad}
d_f(A,C)\leq d_f(A,B)\quad\text{and}\quad d_f(C,B)\leq d_f(A,B)
\end{align}
whenever $A\leq C\leq B$ are Hermitian matrices in $\Mn$ with spectra in $J$.
We denote the class of such functions by $\mathcal R_n(J)$ and write
$\mathcal R_n=\mathcal R_n((0,\infty))$.

  The  classes $\mathcal{M}_n$ and $\mathcal{R}_n$   are nested  as the matrix size increases, i.e.,
\begin{align*}
\mathcal{M}_{n+1}(J)\subseteq\mathcal{M}_n(J)\quad\text{and}\quad \mathcal{R}_{n+1}(J)
\subseteq\mathcal{R}_n(J).
\end{align*}


\section{Exact separation of finite matrix orders}

We construct a family of rational functions whose exact matrix-monotonicity and radial-convexity orders are $m$ and $2m-1$, respectively

\begin{theorem}\label{th-gap}
For every integer $m\ge2$, there exists $\varepsilon_m>0$ such that, for all
$0<\varepsilon<\varepsilon_m$, the rational function
\begin{align}\label{eq--fam}
f_{m,\varepsilon}(t)=\sum_{j=1}^{m}\frac{t+1}{jt+j+2}-\varepsilon\frac{t+1}{(m+1)t+m+3},
\qquad t>0,
\end{align}
satisfies
\begin{align*}
f_{m,\varepsilon}\in\mathcal M_m\setminus\mathcal M_{m+1}\qquad\text{and}\qquad
f_{m,\varepsilon}\in\mathcal R_{2m-1}\setminus\mathcal R_{2m}.
\end{align*}
\end{theorem}


\begin{corollary}\label{cor-gap}
For every integer $n\ge3$,
\begin{align*}
\mathcal{R}_n\not\subseteq\mathcal{M}_n.
\end{align*}
\end{corollary}

\begin{proof}
Fix $n\ge3$ and let
$m=\lceil(n+1)/2\rceil$.
Then $m+1\le n\le2m-1$.
By Theorem~\ref{th-gap}, for sufficiently
small $\varepsilon>0$,
\begin{align*}
f_{m,\varepsilon}
\in\mathcal{R}_{2m-1}
\setminus\mathcal{M}_{m+1}.
\end{align*}
Since both classes are nested as the matrix
size increases, we obtain
$f_{m,\varepsilon}\in
\mathcal{R}_n\setminus\mathcal{M}_n$.
\end{proof}

\subsection{Lower bounds for the two orders}

To prove Theorem~\ref{th-gap}, we use the Cayley type transformation
\begin{align*}
\phi(t)=\frac{t-1}{t+1},\quad t>0,
\end{align*}
which maps $(0,\infty)$ onto $(-1,1)$. For $j\ge1$, put
\begin{align*}
r_j(s)=\frac{1}{j+1-s},\quad s\in(-1,1),
\end{align*}
and define
\begin{align}\label{eq--res}
g_{m,\varepsilon}(s)=\sum_{j=1}^m r_j(s)-\varepsilon r_{m+1}(s).
\end{align}
Then $f_{m,\varepsilon}=g_{m,\varepsilon}\circ\phi$.

Both $\phi$ and its inverse $\phi^{-1}(s)=(1+s)/(1-s)$ are operator
monotone on their   domains. Thus $\phi$ is an order isomorphism between positive definite matrices and Hermitian matrices with spectra in $(-1,1)$. Consequently, for every $n$, $f_{m,\varepsilon}$ is matrix monotone of order $n$ if and only if $g_{m,\varepsilon}$ is. The same holds for Frobenius radial convexity. Thus it suffices to prove the theorem
for $g_{m,\varepsilon}$ on $(-1,1)$.


\begin{lemma}\label{lem-int}
Let $J\subseteq\mathbb R$ be a compact interval,
and let $b_1,\ldots,b_k,h$ be real rational functions with no poles on $J$.
Suppose there is a polynomial $q$, with no zeros on $J$, such that
$qb_1,\ldots,qb_k$ form a basis of the polynomials of degree at most $k-1$.
Then there exists $\kappa>0$ such that
\begin{align}\label{eq--int}
h(H)Xh(H)\le \kappa\sum_{j=1}^k b_j(H)Xb_j(H)
\end{align}
for every Hermitian $H\in\mathbb M_k$ with spectrum in $J$ and every $X\ge0$.
The constant $\kappa$ is independent of $H$ and $X$.
\end{lemma}

\begin{proof}
Put $p_j=qb_j$ and $\widetilde h=qh$.
Let $\lambda_1,\ldots,\lambda_k\in J$ be the
eigenvalues of $H$, counted with multiplicity,
and let $p$ be the Hermite interpolation
polynomial of $\widetilde h$ at these points,
of degree at most $k-1$. Then
\begin{align*}
p(s)=\sum_{\ell=0}^{k-1}
\widetilde h^{[\ell]}
(\lambda_1,\ldots,\lambda_{\ell+1})
\prod_{j=1}^{\ell}(s-\lambda_j),
\end{align*}
see \cite[Section~2.1.5, Eq.~(2.1.5.8)]{SB2002}. By the mean value theorem for divided
differences, there exists $\xi\in J$ such that
\begin{align}\label{eq-mv}
\left|\widetilde h^{[\ell]}(\lambda_1,\ldots,\lambda_{\ell+1})
\right|=\frac{|\widetilde h^{(\ell)}(\xi)|}{\ell!}
\le\frac{1}{\ell!}\max_{s\in J}|\widetilde h^{(\ell)}(s)|.
\end{align}
Since $J$ is compact, the coefficients of the polynomials
$\prod_{j=1}^{\ell}(s-\lambda_j)$ are uniformly bounded.
Together with\eqref{eq-mv}, this shows that the coefficients of $p$
are bounded independently of interpolation points.
Since  $p_1,\ldots,p_k$ form a fixed basis,
we can write
\begin{align*}
p=\sum_{j=1}^k\alpha_jp_j,
\qquad
\sum_{j=1}^k|\alpha_j|^2\le\kappa,
\end{align*}
where $\kappa$ is independent of $H$.

Since $p$ and $\widetilde h$ agree at
the eigenvalues of $H$, the spectral
theorem gives $p(H)=\widetilde h(H)$.
Since $q$ has no zero in $J$, the matrix $q(H)$ is invertible and  we obtain
\begin{align}\label{eq-hbj}
h(H)=\sum_{j=1}^k\alpha_jb_j(H).
\end{align}
Now for any $X\ge0$ and $v\in\mathbb C^k$ we have
\begin{align*}
\langle h(H)Xh(H)v,v\rangle&=\left\|\sum_{j=1}^k\alpha_jX^{1/2}b_j(H)v
\right\|^2\qquad(\text{by \eqref{eq-hbj}})\\
&\le\left(\sum_{j=1}^k|\alpha_j|\,\left\|X^{1/2}b_j(H)v\right\|\right)^2\qquad(\text{by the triangle inequality})\\
&\le \left(\sum_{j=1}^k|\alpha_j|^2\right)\sum_{j=1}^k\left\|X^{1/2}b_j(H)v\right\|^2   \qquad(\text{By the Cauchy--Schwarz inequality})\\
&\le \kappa\sum_{j=1}^k \left\langle b_j(H)Xb_j(H)v,v\right\rangle.
\end{align*}
This proves \eqref{eq--int}.
\end{proof}

\begin{lemma}\label{lem-rsp}
Let $m\ge2$ and put
\begin{align*}
q_m(s)=\prod_{j=1}^m(j+1-s).
\end{align*}
Consider the spaces
\begin{align*}
\mathcal W_m=\spn\{r_1,\ldots,r_m\}\qquad\text{and}\qquad
\mathcal V_m=\spn\{r_ir_j:1\le i,j\le m\}.
\end{align*}
Then
\begin{align}\label{eq--rsp}
q_m\mathcal W_m=\mathcal P_{m-1}\qquad\text{and}\qquad q_m^2\mathcal V_m=\mathcal P_{2m-2},
\end{align}
where $\mathcal P_k$ denotes the space of real polynomials of degree at most $k$.
In particular,
\begin{align*}
\dim\mathcal W_m=m \qquad\text{and}\qquad\dim\mathcal V_m=2m-1.
\end{align*}
\end{lemma}

\begin{proof}
For $1\le j\le m$ we have
\begin{align*}
q_m(s)r_j(s)=\prod_{\substack{1\le i\le m\\i\ne j}}(i+1-s).
\end{align*}
These are polynomials of degree $m-1$, and they are linearly independent, because
at $s=j+1$, all except the $j$th vanish, while the $j$th is nonzero.
 Hence, they form a basis of $\mathcal P_{m-1}$.

Next, note that for $i\ne j$ we can write
\begin{align*}
r_ir_j=\frac{(j-1)r_1r_j-(i-1)r_1r_i}{j-i}.
\end{align*}
Thus $\mathcal V_m$ is spanned by
\begin{align*}
r_1^2,\ldots,r_m^2,\quad r_1r_2,\ldots,r_1r_m.
\end{align*}
Moreover, these $2m-1$ functions are linearly
independent.
Suppose that
\begin{align}\label{eq-sajbj}
\sum_{j=1}^m a_jr_j^2+\sum_{j=2}^m b_jr_1r_j=0.
\end{align}
For each $j=1,\ldots,m$, we have
\begin{align*}
a_j=\lim_{s\to j+1}(j+1-s)^2\left(\sum_{i=1}^m a_ir_i(s)^2
+\sum_{i=2}^m b_ir_1(s)r_i(s)\right)=0.
\end{align*}
Together with \eqref{eq-sajbj} this gives
$$\sum_{j=2}^m b_jr_1r_j=0.$$
 For $2\le j\le m$, multiplying by $(2-s)(j+1-s)$ and letting $s\to j+1$
gives $b_j=0$. Hence the functions are linearly independent.
Consequently, we get $\dim\mathcal V_m=2m-1$.

Each $q_m^2r_ir_j$ is a polynomial of degree
at most $2m-2$. Hence
\begin{align*}
q_m^2\mathcal V_m\subseteq\mathcal P_{2m-2}.
\end{align*}
Since $q_m$ is a nonzero polynomial,
multiplication by $q_m^2$ preserves dimension.
Thus
\begin{align*}
\dim(q_m^2\mathcal V_m)=2m-1=\dim\mathcal P_{2m-2},
\end{align*}
and  so $q_m^2\mathcal V_m=\mathcal P_{2m-2}$.
\end{proof}


\begin{lemma}\label{lem-dom}
For every $m\ge2$, there exist constants
$\kappa_m^{(1)},\kappa_m^{(2)}>0$, depending only on $m$, such that the following
inequalities hold for all Hermitian matrices $H$ with spectrum in $(-1,1)$ and every $X\ge0$.
\begin{enumerate}
\item[(i)] If $H\in\mathbb M_m$, then
\begin{align}\label{eq--dom1}
r_{m+1}(H)Xr_{m+1}(H)\le \kappa_m^{(1)}\sum_{j=1}^m r_j(H)Xr_j(H).
\end{align}

\item[(ii)] If $H\in\mathbb M_{2m-1}$, then
\begin{align}\label{eq--dom2}
\sum_{i=1}^mr_i(H)r_{m+1}(H)Xr_i(H)r_{m+1}(H)\le \kappa_m^{(2)}
\sum_{i,j=1}^mr_i(H)r_j(H)Xr_i(H)r_j(H).
\end{align}
\end{enumerate}
\end{lemma}

\begin{proof}
We apply Lemma~\ref{lem-int} on $[-1,1]$.
Recall that
\begin{align*}
q_m(s)=\prod_{j=1}^m(j+1-s),
\end{align*}
which has no zeros on this interval.  By Lemma~\ref{lem-rsp}, the polynomials
$q_mr_1,\ldots,q_mr_m$ form a basis of $\mathcal P_{m-1}$. Applying
Lemma~\ref{lem-int} with $q=q_m$, $h=r_{m+1}$ and $b_j=r_j$ proves (i).

For (ii), put $k=2m-1$ and let $b_1,\ldots,b_k$ be the basis of
$\mathcal V_m$ given by
\begin{align*}
r_1^2,\ldots,r_m^2, \quad r_1r_2,\ldots,r_1r_m.
\end{align*}
By Lemma~\ref{lem-rsp}, $q_m^2b_1,\ldots,q_m^2b_k$ form a basis
of $\mathcal P_{k-1}$.

For each $1\le i\le m$, apply Lemma~\ref{lem-int} with $q=q_m^2$
and $h=r_ir_{m+1}$ to find $\kappa_i>0$ such that
\begin{align}\label{eq-rirm+1}
r_i(H)r_{m+1}(H)Xr_i(H)r_{m+1}(H)\le\kappa_i\sum_{\ell=1}^kb_\ell(H)Xb_\ell(H).
\end{align}
Since $X\ge0$ and each $b_\ell$ is one of
the products $r_ir_j$, we have
\begin{align}\label{eq-blH}
\sum_{\ell=1}^k b_\ell(H)Xb_\ell(H)\le\sum_{i,j=1}^m
r_i(H)r_j(H)Xr_i(H)r_j(H).
\end{align}
Now, (ii) follows from \eqref{eq-rirm+1}, \eqref{eq-blH}, and summing over $i=1,\ldots,m$, with
\begin{align*}
\kappa_m^{(2)}=\sum_{i=1}^m\kappa_i.
\end{align*}
\end{proof}


\begin{proposition}\label{prop-mon}
Let $m\ge2$. If $0<\varepsilon\le(\kappa_m^{(1)})^{-1}$,
then the function $g_{m,\varepsilon}$ is matrix monotone of order $m$ on $(-1,1)$.
\end{proposition}

\begin{proof}
Let $H\in\mathbb M_m$ be Hermitian with spectrum in $(-1,1)$, and let $X\ge0$.
Since
\begin{align*}
\mathrm{D}r_j(H)[X]=r_j(H)Xr_j(H),
\end{align*}
Lemma~\ref{lem-dom}(i) gives
\begin{align}\label{eq-devr}
\mathrm{D}g_{m,\varepsilon}(H)[X]&=\sum_{j=1}^m r_j(H)Xr_j(H)
-\varepsilon r_{m+1}(H)Xr_{m+1}(H)\\
&\ge(1-\varepsilon\kappa_m^{(1)})\sum_{j=1}^m r_j(H)Xr_j(H)\ge0\nonumber
\end{align}
in which    positivity follows from our assumption on $\varepsilon$.

Now let $A\le B$ be Hermitian matrices
with spectra in $(-1,1)$, and put
$H_t=A+t(B-A)$ for $0\le t\le1$.
Then $H_t$ also has spectrum in $(-1,1)$,
and hence
\begin{align*}
g_{m,\varepsilon}(B)-g_{m,\varepsilon}(A)
&=\int_0^1\mathrm{D}g_{m,\varepsilon}(H_t)[B-A]\,\mathrm{d}t\ge0.
\end{align*}
Thus $g_{m,\varepsilon}$ is matrix monotone
of order $m$.
\end{proof}


\begin{lemma}\label{lem-inc}
Let $A,B\in\mathbb M_n$ be Hermitian with spectra in $(-1,1)$, and suppose
that $A\le B$. Put $X=B-A$. Then
\begin{align}\label{eq--inc}
\frac13 r_j(E)Xr_j(E)\le r_j(B)-r_j(A)\le 3\,r_j(E)Xr_j(E)
\end{align}
  for every $j\ge1$ and $E\in\{A,B\}$.
\end{lemma}

\begin{proof}
Put $\Delta_j=r_j(B)-r_j(A)$.
The resolvent identity gives
\begin{align*}
\Delta_j=r_j(A)Xr_j(B)=r_j(B)Xr_j(A).
\end{align*}
Consequently,
\begin{align}\label{eq-dl1}
\Delta_j=r_j(A)(X+Xr_j(B)X)r_j(A)=r_j(B)(X-Xr_j(A)X)r_j(B).
\end{align}
Since $0\le X\le2I$ and $r_j(B)\le j^{-1}I$, we have
\begin{align}\label{eq-xrb}
0\le Xr_j(B)X\le\frac1j X^2\le2X.
\end{align}
It follows from \eqref{eq-xrb} and  the first equality of \eqref{eq-dl1} that
\begin{align*}
r_j(A)Xr_j(A)\le\Delta_j\le3\,r_j(A)Xr_j(A),
\end{align*}
proving the right side inequality in \eqref{eq--inc}.
Since $A\ge-I$, we have $r_j(A)\ge(j+2)^{-1}I$ and so
\begin{align}\label{eq-ner1}
r_j(A)(I-A)=I-jr_j(A)\le\frac{2}{j+2}I.
\end{align}
Moreover,  $X\le I-A$ implies that
\begin{align}\label{eq-ner2}
r_j(A)^{1/2}Xr_j(A)^{1/2}&\le r_j(A)(I-A).
\end{align}
From \eqref{eq-ner1} and \eqref{eq-ner2} we get
\begin{align*}
\|r_j(A)^{1/2}X^{1/2}\|^2\le\frac{2}{j+2},
\end{align*}
and hence
\begin{align*}
0\le Xr_j(A)X\le\frac{2}{j+2}X.
\end{align*}
The second  equality of \eqref{eq-dl1} now yields
\begin{align*}
\frac{j}{j+2}r_j(B)Xr_j(B)\le\Delta_j\le r_j(B)Xr_j(B).
\end{align*}
As $j/(j+2)\ge1/3$, this proves the left side  inequality in \eqref{eq--inc}.
\end{proof}


\begin{proposition}\label{prop-rad}
Let $m\ge2$. If
$0<\varepsilon\le(18\kappa_m^{(2)})^{-1}$,
then $g_{m,\varepsilon}$ is Frobenius
radially convex of order $2m-1$ on $(-1,1)$.
\end{proposition}

\begin{proof}
Put $n=2m-1$ and $g=g_{m,\varepsilon}$.
Let $A\le B$ be Hermitian matrices in $\mathbb M_n$ with spectra in $(-1,1)$,
and set
\begin{align*}
X=B-A,\qquad\Delta_j=r_j(B)-r_j(A).
\end{align*}
Fix $E\in\{A,B\}$ and write $R_j=r_j(E)$. By Lemma~\ref{lem-inc} we have
\begin{align}\label{eq-lmin}
\frac13 R_jXR_j\le\Delta_j\le3R_jXR_j.
\end{align}

Since
\begin{align*}
g(B)-g(A)
=\sum_{i=1}^m\Delta_i-\varepsilon\Delta_{m+1},
\end{align*}
the   derivative formula \eqref{eq-devr} gives
\begin{align*}
\mathrm Dg(E)[g(B)-g(A)]&=\sum_{j=1}^m r_j(E)(g(B)-g(A))r_j(E)
-\varepsilon r_{m+1}(E)(g(B)-g(A))r_{m+1}(E)\\
&=\sum_{i,j=1}^m R_j\Delta_iR_j-\varepsilon\sum_{i=1}^m
\bigl(R_i\Delta_{m+1}R_i
+R_{m+1}\Delta_iR_{m+1}\bigr)+\varepsilon^2
R_{m+1}\Delta_{m+1}R_{m+1}.
\end{align*}
The last term is positive semidefinite. Using \eqref{eq-lmin} for the
remaining terms, we obtain
\begin{align}\label{eq-derg}
\begin{split}
\mathrm Dg(E)[g(B)-g(A)]&\ge\frac13\sum_{i,j=1}^m
R_iR_jXR_iR_j-6\varepsilon\sum_{i=1}^mR_iR_{m+1}XR_iR_{m+1}\\
&\ge \left(\frac13-6\varepsilon\kappa_m^{(2)}\right)
\sum_{i,j=1}^m R_iR_jXR_iR_j
\ge0,
\end{split}
\end{align}
where the second inequality follows from Lemma~\ref{lem-dom}(ii).

Now let $A\le C\le B$ and  put
$B_t=C+t(B-C)$ for $0\le t\le1$. Then $A\le B_t$  and \eqref{eq-derg} gives
\begin{align}\label{eq-pdg}
\mathrm Dg(B_t)[g(B_t)-g(A)]\ge0.
\end{align}
Since $\mathrm Dg(B_t)$ is self-adjoint, we have
\begin{align*}
 \frac{\mathrm d}{\mathrm dt}\|g(B_t)-g(A)\|_F^2
  =2\tr\bigl((B-C)\mathrm Dg(B_t)[g(B_t)-g(A)]\bigr).
\end{align*}
From  $B-C\ge0$ and \eqref{eq-pdg}, we infer that the above
  derivative is nonnegative, and so  the squared Frobenius norm is
nondecreasing in $t$. Hence
\begin{align*}
\|g(C)-g(A)\|_F\le\|g(B)-g(A)\|_F.
\end{align*}

Furthermore, considering $A_t=A+t(C-A)$ and applying the same
argument  to $A_t\le B$ at the endpoint $A_t$, we obtain
\begin{align*}
\frac{\mathrm d}{\mathrm dt}\|g(B)-g(A_t)\|_F^2=-2\operatorname{Tr}\bigl(
(C-A)\mathrm Dg(A_t)[g(B)-g(A_t)]\bigr),
\end{align*}
which is nonpositive, since $C-A\geq0$ and $\mathrm{D}g(A_t)[g(B)-g(A_t)]\geq0$. Consequently,
\begin{align*}
\|g(B)-g(C)\|_F
\le\|g(B)-g(A)\|_F.
\end{align*}
Thus $g$ is Frobenius radially convex
of order $2m-1$.
\end{proof}



\subsection{Sharpness of the bounds}

In this subsection, we prove that $g_{m,\varepsilon}$ fails
matrix monotonicity at order $m+1$ for every $\varepsilon>0$, and Frobenius radial convexity at order $2m$ for sufficiently small $\varepsilon>0$.

\begin{proposition}\label{prop-nmon}
For every $m\ge2$ and $\varepsilon>0$,
the function $g_{m,\varepsilon}$ is not
matrix monotone of order $m+1$ on $(-1,1)$.
\end{proposition}

\begin{proof}
Choose distinct points
$s_1,\ldots,s_{m+1}\in(-1,1)$. Since
\begin{align}\label{eq-fddr}
r_j^{[1]}(s,t)=r_j(s)r_j(t),
\end{align}
the Loewner matrix of $g_{m,\varepsilon}$
at these points is
\begin{align}\label{eq-lng}
\begin{split}
L_{g_{m,\varepsilon}}(s_1,\ldots,s_{m+1})&=\sum_{j=1}^mL_{r_j}(s_1,\ldots,s_{m+1})
-\varepsilon L_{r_{m+1}}(s_1,\ldots,s_{m+1})\\
&=\left[
\sum_{j=1}^m r_j(s_k)r_j(s_\ell)
-\varepsilon r_{m+1}(s_k)r_{m+1}(s_\ell)
\right]_{k,\ell=1}^{m+1}\qquad(\text{by \eqref{eq-fddr}})\\
&=U\operatorname{diag} (1,\ldots,1,-\varepsilon)U^*,
\end{split}
\end{align}
where
\begin{align*}
U=\bigl[r_j(s_i)\bigr]_{i,j=1}^{m+1}.
\end{align*}

We show that $U$ is invertible. Suppose that $Uc=0$ for some $c\in\mathbb R^{m+1}$.
Define
\begin{align*}
p(s)=q_{m+1}(s)\sum_{j=1}^{m+1}c_jr_j(s)
\end{align*}
in which $q_{m+1}$ is defined as in Lemma~\ref{lem-rsp}.
Since $p$ has degree at most $m$ and vanishes at the $m+1$ distinct points
$s_1,\ldots,s_{m+1}$, it follows that $p=0$.

By Lemma~\ref{lem-rsp}, the polynomials $q_{m+1}r_1,\ldots,q_{m+1}r_{m+1}$
are linearly independent, so $c=0$.

Thus $U$ is invertible. Using  Sylvester's law of inertia and noting \eqref{eq-lng},  the Loewner matrix of $g_{m,\varepsilon}$ has exactly one negative eigenvalue.
Loewner's criterion therefore shows that $g_{m,\varepsilon}$ is not matrix
monotone of order $m+1$.
\end{proof}


The following necessary condition for Frobenius
radial convexity follows from
\cite[Propositions 4.5 and 4.8]{CFM2026}, whose proofs extend to arbitrary open intervals.

\begin{lemma}\label{lem-had}
Let $J\subseteq\mathbb R$ be an open interval
and $f\in C^1(J)$. If $f\in\mathcal R_n(J)$,
then
\begin{align}
\label{eq--had}
L_f(s_1,\ldots,s_n)\circ
L_f(s_1,\ldots,s_n)\ge0
\end{align}
for all $s_1,\ldots,s_n\in J$.
\end{lemma}


\begin{proposition}\label{prop-nrad}
For every $m\ge2$, the function
$g_{m,\varepsilon}$ is not Frobenius radially
convex of order $2m$ on $(-1,1)$ for all
sufficiently small $\varepsilon>0$.
\end{proposition}

\begin{proof}
Fix $m\ge2$, put $n=2m$, and choose distinct
points $s_1,\ldots,s_n\in(-1,1)$. Set
\begin{align*}
p(s)=\prod_{\ell=1}^n(s-s_\ell)\qquad\text{and}\qquad
z_\ell=\frac{q_m(s_\ell)^2}{p'(s_\ell)},\qquad 1\le\ell\le n,
\end{align*}
and write $z=(z_1,\ldots,z_n)^T$. For
$1\le i,j\le m+1$, let
\begin{align*}
u_j=(r_j(s_1),\ldots,r_j(s_n))^T\qquad\text{and}\qquad w_{ij}=u_i\circ u_j.
\end{align*}

By the Lagrange interpolation formula
\cite[Section 2.1.1]{SB2002}, every polynomial
$h$ of degree at most $n-1$ satisfies
\begin{align*}
h(s)=\sum_{\ell=1}^n
h(s_\ell)\frac{p(s)}
{(s-s_\ell)p'(s_\ell)}.
\end{align*}
Now if  $ \deg h\le n-2$, then comparing the coefficients of $s^{n-1}$  in both sides of above identity   gives
\begin{align}\label{eq-cosn-1}
\sum_{\ell=1}^n \frac{h(s_\ell)}{p'(s_\ell)}=0.
\end{align}

For $1\le i,j\le m$, the function $q_m(s)^2r_i(s)r_j(s)$  is a polynomial of degree $n-2$. Hence, by \eqref{eq-cosn-1} we have
\begin{align}\label{eq--orth}
z^*w_{ij}&=\sum_{\ell=1}^n\frac{q_m(s_\ell)^2
r_i(s_\ell)r_j(s_\ell)}{p'(s_\ell)}=0,\qquad 1\le i,j\le m.
\end{align}

For $1\le i\le m$, the function
$q_m^2r_i$ is a polynomial of degree
$n-1$. Evaluating its Lagrange interpolation
formula at $s=m+2$ gives
\begin{align*}
q_m(m+2)^2r_i(m+2)
=p(m+2)\sum_{\ell=1}^n
\frac{q_m(s_\ell)^2r_i(s_\ell)}
{(m+2-s_\ell)p'(s_\ell)}.
\end{align*}
Since $p(m+2)\ne0$, we obtain
\begin{align}\label{eq--mix}
z^*w_{i,m+1}=\sum_{\ell=1}^n\frac{q_m(s_\ell)^2r_i(s_\ell)}
{(m+2-s_\ell)p'(s_\ell)}=\frac{q_m(m+2)^2r_i(m+2)}{p(m+2)}\ne0.
\end{align}

Let $L_\varepsilon$ denote the Loewner matrix
of $g_{m,\varepsilon}$ at $s_1,\ldots,s_n$.
By \eqref{eq-fddr}, we have
\begin{align*}
L_\varepsilon=\sum_{j=1}^m u_ju_j^*-\varepsilon u_{m+1}u_{m+1}^*.
\end{align*}
Noticing $(u_i u_i^*)\circ(u_j u_j^*)=w_{ij}w_{ij}^*$, we infer
\begin{align*}
L_\varepsilon\circ L_\varepsilon=\sum_{i,j=1}^m w_{ij}w_{ij}^*
-2\varepsilon\sum_{i=1}^mw_{i,m+1}w_{i,m+1}^*+\varepsilon^2
w_{m+1,m+1}w_{m+1,m+1}^*.
\end{align*}
Using \eqref{eq--orth}, we obtain
\begin{align*}
z^*(L_\varepsilon\circ L_\varepsilon)z=-2\varepsilon
\sum_{i=1}^m|z^*w_{i,m+1}|^2+\varepsilon^2|z^*w_{m+1,m+1}|^2.
\end{align*}
By \eqref{eq--mix}, the coefficient of $-2\varepsilon$ is strictly positive.
Thus $z^*(L_\varepsilon\circ L_\varepsilon)z<0$, i.e., $L_\varepsilon\circ L_\varepsilon$
is not positive semidefinite for all sufficiently small
$\varepsilon>0$.  By Lemma~\ref{lem-had}, the function $g_{m,\varepsilon}$ is not Frobenius
radially convex of order $2m$.
\end{proof}


\begin{proof}[Proof of Theorem~\ref{th-gap}]
Fix $m\ge2$. By Proposition~\ref{prop-nrad},
there exists $\delta_m>0$ such that
$g_{m,\varepsilon}$ is not Frobenius radially
convex of order $2m$ whenever
$0<\varepsilon<\delta_m$.
Set
\begin{align*}
\varepsilon_m
=\min\left\{
\frac{1}{\kappa_m^{(1)}},
\frac{1}{18\kappa_m^{(2)}},
\delta_m
\right\}.
\end{align*}
For every $0<\varepsilon<\varepsilon_m$,
Propositions~\ref{prop-mon}, \ref{prop-rad},
\ref{prop-nmon}, and \ref{prop-nrad} give
\begin{align*}
g_{m,\varepsilon}
&\in\mathcal M_m((-1,1))
\setminus\mathcal M_{m+1}((-1,1)),\\
g_{m,\varepsilon}
&\in\mathcal R_{2m-1}((-1,1))
\setminus\mathcal R_{2m}((-1,1)).
\end{align*}
Since $f_{m,\varepsilon}
=g_{m,\varepsilon}\circ\phi$,
the order-preserving change of variables
established above gives the corresponding
assertions for $f_{m,\varepsilon}$ on
$(0,\infty)$.
\end{proof}

\begin{remark}
The counterexamples in Corollary~\ref{cor-gap}
are not  matrix monotone decreasing  of any order.
 Indeed, for
$0<\varepsilon<\varepsilon_m$, applying
\eqref{eq-devr} with $H=sI_m$ and $X=I_m$
gives
$$g_{m,\varepsilon}'(s)\ge(1-\varepsilon\kappa_m^{(1)})
\sum_{j=1}^m r_j(s)^2>0,\qquad -1<s<1.$$
Since $\phi'(t)>0$ for $t>0$, the function
$f_{m,\varepsilon}=g_{m,\varepsilon}\circ\phi$
is strictly increasing on $(0,\infty)$.
Consequently, the functions used to show $\mathcal R_n\not\subseteq\mathcal M_n$
in Corollary~\ref{cor-gap} are neither
matrix monotone increasing nor matrix
monotone decreasing of order $n$.
\end{remark}

\bibliographystyle{amsplain}

\end{document}